\documentclass[12pt]{article}
\usepackage{theorem}
\usepackage{graphics}
\newtheorem{lemma}{Lemma}

\newtheorem{theorem}[lemma]{Theorem}

{\theorembodyfont{\upshape}}
{\theorembodyfont{\upshape}}
{\theorembodyfont{\upshape}}
{\theorembodyfont{\upshape}}
{\theorembodyfont{\upshape}}
{\theorembodyfont{\upshape}}
{\theorembodyfont{\upshape}}

\newcommand{\R}{{\bf R}}

\newcommand{\rme}{{\rm e}}
\newcommand{\rmd}{\,{\rm d}}

\newcommand{\sig}{\sigma}
\newcommand{\alp}{\alpha}
\newcommand{\bet}{\beta}
\newcommand{\gam}{\gamma}
\newcommand{\lam}{\lambda}
\newcommand{\del}{\delta}

\newcommand{\Gam}{\Gamma}
\newcommand{\kap}{\kappa}

\newcommand{\norm}{\Vert}
\renewcommand{\Re}{{\rm Re}\,}

\newcommand{\Proof}{\underbar{Proof}{\hskip 0.1in}}
\newcommand{\Note}{\underbar{Note}{\hskip 0.1in}}

\newcommand{\Schrodinger}{Schr\"odinger }

\newcommand{\la}{{\langle}}
\newcommand{\ra}{{\rangle}}
\newcommand{\pr}{\prime}

\newcommand{\nopic}[1]{}

\newcommand{\be}{\begin{equation}}
\newcommand{\ee}{\end{equation}}

\usepackage{hyperref}
\title{ASYMPTOTIC BEHAVIOUR OF\\
QUASI-ORTHOGONAL POLYNOMIALS}
\author{E.B. Davies}
\date{8 April 2003}
\begin{document}
\maketitle

\begin{abstract}
We obtain explicit upper and lower bounds on the norms of the
spectral projections of the non-self-adjoint harmonic oscillator.
Some of our results apply to a variety of other families of
orthogonal polynomials.
\end{abstract}

\section{Introduction}

We consider polynomials $p_n$ which are orthogonal with respect to
a complex weight $\sig$ on $[0,\infty)$ in the following sense. We
suppose that $p_n$ is of degree $n$ and
\[
\int_0^\infty p_m(x)p_n(x)\sig(x)^2\,\rmd x=\del_{m,n}
\]
for all non-negative integers $m,\, n$. (All of our statements and
proofs can be rewritten with $(0,\infty)$ replaced by $\R$, and we
will not keep repeating this point.) If $\sig>0$ and $p_m$ are
real-valued, then they are orthonormal in $L^2((0,\infty),
\sig(x)^2\,\rmd x)$ in the usual sense, but for complex-valued
$\sig$ such an interpretation is not possible. Our goal is to
obtain bounds on the quantities
\[
N_n=\int_0^\infty |p_n(x)\sig(x)|^2\,\rmd x
\]
for all $n$.

This problem arose in the context of the non-self-adjoint harmonic
oscillator
\be
(Hf)(x)=-f^{\pr\pr}(x)+z^4x^2f(x)\label{harmonic}
\ee
acting in $L^2(\R)$ for some complex $z$. In this situation the
relevant weight is
\[
\sig(x)=\rme^{-z^2x^2/2}
\]
and $N_n$ is the norm of the spectral projection $P_n$ of $H$
associated with its $n$th eigenvalue, $\lam_n=z^{2}(2n+1)$. In the
numerical literature $N_n$ is called the condition number of the
eigenvalue $\lam_n$. Numerical calculations in \cite{AD} indicated
that $\norm P_n\norm$ increases at an exponential rate as
$n\to\infty$, and it was proved in \cite{EBD3} that there was no
polynomial bound on $\norm P_n\norm$ for this and certain other
\Schrodinger operators. The super-polynomial rate of increase of
the associated resolvent norms in the semi-classical limit was
proved in \cite{EBD2} by a method which was greatly generalized in
\cite{Zw}. For certain classes of operators with analytic
coefficients it was recently proved that the resolvent norms
increase at an exponential rate in the semiclassical limit,
\cite{DSZ}. However, the precise exponential constants have not
been identified in any example.

A consequence of our theorems is that there exists a positive
critical constant $t_z$ such that the `spectral expansion'
\[
\rme^{-Ht}=\sum_{n=0}^\infty \rme^{-\lam_nt}P_n
\]
is norm convergent if $t>t_z$ and divergent if $0\leq t<t_z$. Our
method provides explicit upper and lower bounds on $t_z$ but not
its precise value.

The problem may be reformulated as finding the norms of
$\phi_n(x)=p_n(x)\sig(x)$ in $L^2((0,\infty),\rmd x)$, where
$\phi_n$ are obtained by applying a modified Gram-Schmidt
orthogonalization process to the functions $x^n\sig(x)$. This
procedure is modified in the sense that we require
\[
\int_0^\infty \phi_m(x)\phi_n(x)\,\rmd x=\del_{m,n}
\]
without any complex conjugates. This is equivalent to requiring
that $\phi_m$ and $\phi_n^\ast(x)=\overline{\phi_n(x)}$ form a
biorthogonal system in $L^2((0,\infty),\rmd x)$ in the sense that
\[
\la \phi_m,\phi_n^\ast\ra=\del_{m,n}
\]
for all $m,\,n$. If $P_n$ is the (non-orthogonal) projection
\[
P_nf=\la f,\phi_n^\ast\ra \phi_n
\]
then $P_mP_n=\del_{m,n}P_n$ for all $m,\,n$ and it is easily seen
that
\[
\norm P_n \norm =N_n.
\]

In order to make some progress with this problem, we make the
following assumptions on the weight $\sig$. We assume that
$\sig(z)$ is an analytic function of $z$ in the sector $S=\{
z:|\arg(z)|<\alp\}$, and that it is positive on the real axis. We
also assume that
\[
\int_0^\infty x^n |\sig(\rme^{i\theta}x)|^2\,\rmd x<\infty
\]
for all $n\geq 0$ and $|\theta|<\alp$, in order that $p_n$ should
be well-defined. Our most important condition is that
\be
|\sig(\rme^{i\theta} r)|\geq c_\theta\sig(s_\theta r)\label{basic}
\ee
for  all $|\theta|<\alp$ and all $r>0$, where $c_\theta>0$ and $0<
s_\theta <1$. Our main theorem provides a lower bound on $N_n$ for
the weight $x\to\sig(\rme^{i\theta}x)$ under these assumptions.
Examples of such weights are given in Section 3. Finally in
Section 4 we compare the bounds obtained with numerical evidence.

\section{The Lower Bound}

Let $\{p_n\}_{n=0}^\infty$ denote the standard orthonormal
sequence of real-valued polynomials with respect to the positive
weight $\sig^2$ on $(0,\infty)$. We define
\[
p_{n,z}(x)=z^{1/2}p_n(zx)
\]
where $z\in S$ and $x>0$. If $z>0$ then
\[
\int_0^\infty p_{m,z}(x)p_{n,z}(x)\sig(zx)^2\,\rmd x=\del_{m,n}
\]
by making the change of variable $zx=u$. By analytic continuation
the same holds for all complex $z\in S$. We are interested in
obtaining a lower bound on the quantity
\[
N_{n,z}=\int_0^\infty |p_{n,z}(x)\sig(zx)|^2\,\rmd x
\]
for complex $z\in S$. Note that $N_{n,z}=1$ for all positive real
$z$.

\begin{theorem}\label{main}
Under the assumption (\ref{basic}) we have
\be
N_{n,z}\geq c_\theta^2s_\theta^{-2n-1}\label{lowerbound}
\ee
provided $z=r\rme^{i\theta}$ and $|\theta|<\alp$.
\end{theorem}

\Proof We have
\begin{eqnarray*}
N_{n,z}&= &|z|\int_0^\infty |p_n(zx)\sig(zx)|^2\,\rmd x\\
&\geq & c_\theta^2r\int_0^\infty |p_n(zx)\sig(s_\theta rx)|^2\,\rmd x\\
&=& c_\theta^2s_\theta^{-1}\int_0^\infty |p_n(zx/s_\theta r)\sig(
x)|^2\,\rmd x.
\end{eqnarray*}

Now
\[
p_n(zx/s_\theta r)
=z^ns_\theta^{-n}r^{-n}p_n(x)+\sum_{j=0}^{n-1}k_jp_j(x)
\]
for constants $k_j$ which we need not evaluate. By the
orthogonality of the polynomials, we have
\begin{eqnarray*}
\int_0^\infty |p_n(zx/s_\theta r)\sig( x)|^2\,\rmd x
&=&s_\theta^{-2n} +\sum_{j=0}^{n-1} |k_j|^2\\
 &\geq& s_\theta^{-2n}.
\end{eqnarray*}
The statement of the theorem follows.

We next consider the example
\[
\sig(z)=z^{\gam/2}\rme^{-z^\bet}
\]
where $\gam>-1$ and $\bet >0$. If $r>0$ and $|\theta|<\pi/(2\bet)$
then
\be
|\sig(r\rme^{i\theta})|=
r^{\gam/2}\rme^{-r^\beta\cos(\theta\beta)}=c_\theta \sig(s_\theta
r)\label{weightidentity}
\ee
where $s_\theta=\{\cos(\theta\bet)\}^{1/\bet}$ and
$c_\theta=s_\theta^{-\gam/2}$. After replacing $(0,\infty)$ by
$(-\infty,\infty)$, the particular choice $\gam=0$ and $\bet=2$
leads one to the study of the Hermite polynomials with a complex
scaling, which is relevant to the non-self-adjoint harmonic
oscillator. The choice $\bet =1$ leads to the Laguerre polynomials
$L_n^{\gam}$. As far as we know, all other choices lead to
non-classical polynomials.

The following theorem provides a more general type of weight
satisfying (\ref{basic}), and can itself easily be generalized.

\begin{theorem} If
\[
\sig(x)=\exp\{-\sum_{j=1}^n c_j x^j\}
\]
for all $x\in (0,\infty)$, where $c_j\in\R$ for all $j$ and
$c_n>0$, then $\sig$ satisfies (\ref{basic}) provided
$|\theta|<\pi/(2n)$.
\end{theorem}

\Proof We have to find $k_\theta >0$ and $s_\theta\in (0,1)$ such
that
\[
\sum_{j=1}^n c_j \cos(j\theta)r^j \leq k_\theta + \sum_{j=1}^n c_j
s_\theta^j r^j
\]
for all $r>0$ and $|\theta|<\pi/(2n)$. The validity of such an
inequality depends upon the coefficient of $r^n$. We achieve the
required bound $ \cos(n\theta)<s_\theta^n <1$  by putting
\[
s_\theta=\left\{ (1+\cos(n\theta))/2  \right\}^{1/n}.
\]

\Note If $(0,\infty)$ is replaced by $\R$ in the above theorem, we
must also assume that $n$ is even.

\section{The Upper Bound}

It is surprisingly difficult to obtain an upper bound on $N_n$,
and we treat only two cases. We start with the orthonormal
sequence of Laguerre polynomials, associated with the weight
$\sig(x)=\rme^{-x/2}$ on $(0,\infty)$. We have
\begin{eqnarray*}
p_n(x)&=&\frac{(-1)^n}{n!}\, \rme^x \frac{\rmd^n}{\rmd x^n} \left(
x^n\rme^{-x}\right)\\
&=& \sum_{r=0}^n b_{n,r}x^r
\end{eqnarray*}
where
\[
b_{n,r}=(-1)^{n-r} \frac{n!}{(r!)^2 (n-r)!}
\]
satisfies
\[
|b_{n,r}|\leq 2^n/r!
\]
by virtue of the general inequality
\[
(r+s)!\leq 2^{r+s}r!\,s!
\]
The following theorem provides an upper bound on $N_{n,z}$ which
complements the lower bound of Theorem~\ref{main}.

\begin{theorem} \label{lagtheo}
If $\sig(x)=\rme^{-x/2}$ and $z=\rme^{i\theta}$
then
\[
N_{n,z}\leq  s_\theta^{-2n-1} 2^{4n+2}
\]
for all $n\geq 0$, provided $|\theta|<\pi/2$ and
$s_\theta=\cos(\theta)$.
\end{theorem}

\Proof We start with the equality
\[
N_{n,z}=s_\theta^{-1} \int_0^\infty
|p_n(\rme^{i\theta}x/s_\theta)\sig(x)|^2\,\rmd x,
\]
which is proved as in Theorem~\ref{main}. We have $c_\theta=1$ and
$s_\theta=\cos(\theta)$ by (\ref{weightidentity}). We deduce that
\begin{eqnarray*}
N_{n,z}&\leq & s_\theta^{-1} \int_0^\infty\sum_{r,s=0}^n
 |b_{n,r}b_{n,s}|s_\theta^{-r-s}
x^{r+s}\rme^{-x}\,\rmd x \\
&\leq &  s_\theta^{-2n-1}2^{2n}\int_0^\infty\sum_{r,s=0}^n
\frac{x^{r+s}}
{ r!s!}\rme^{-x}\,\rmd x \\
&\leq &  s_\theta^{-2n-1}2^{2n}
\sum_{r,s=0}^n \frac{(r+s)!}{r!\, s!}\\
&\leq &  s_\theta^{-2n-1}2^{2n}
\left(\sum_{r=0}^n 2^r\right)^2\\
&\leq & s_\theta^{-2n-1} 2^{4n+2}.
\end{eqnarray*}

\Note This proof can be extended to more general weights provided
suitable bounds on the coefficients $b_{n,r}$ can be obtained, but
in general this is not easy.

We next consider the non-self-adjoint harmonic oscillator. The
orthonormal sequence of polynomials corresponding to the weight
$\sig(x)=\rme^{-x^2/2}$ is given by $p_n(x)=k_n H_n(x)$, where
\[
k_n=\pi^{-1/4}2^{-n/2}(n!)^{-1/2}
\]
and $H_n$ are the Hermite polynomials
\[
H_n(x)=(2x)^n-\frac{n!}{1!\,(n-2)!} (2x)^{n-2}+
\frac{n!}{2!\,(n-4)!} (2x)^{n-4}-...
\]
We will need the following lemma.

\begin{lemma} If $r,\, s$ are non-negative integers then
\[
\int_{-\infty}^\infty x^{2r+2s}\rme^{-x^2}\,\rmd x \leq \pi^{1/2}
2^{r+s}r!\, s!
\]
\end{lemma}

\Proof The left hand-side equals
\begin{eqnarray*}
\int_0^\infty u^{r+s}\rme^{-s}s^{-1/2}\,\rmd s&=&
\Gam(r+s+1/2)\\
&\leq & \pi^{1/2}\Gam(r+s+1)\\
&\leq & \pi^{1/2}2^{r+s} r!\,s!
\end{eqnarray*}

In the following theorem we restrict attention to the case of even
integers; the treatment of odd integers is very similar.

\begin{theorem}\label{upper}
Let $z=\rme^{i\theta}$ where $|\theta|<\pi/4$, and put
$s_\theta=(\cos(2\theta))^{1/2}$. Then
\[
N_{2n,z}\leq \pi(n+1)^{1/2} 2^{4n+2}s_\theta^{-4n-1}.
\]
for all non-negative integers $n$.
\end{theorem}

\Proof We start with the identity
\[
p_{2n}(x)=\sum_{r=0}^n b_{n,r}x^{2r}
\]
where
\[
b_{n,r}=\frac{  (-1)^{n-r}2^{2r-n}\sqrt{(2n)!}  } {
\pi^{1/4}(n-r)!\, (2r)! }.
\]
In the following chain of inequalities we will use
\[
2^{-2r}\sqrt{r+1}\leq \frac{(r!)^2}{(2r)!}\leq
2^{-2r}\sqrt{\pi(r+1)}
\]
for all non-negative integers $r$; this is proved using induction
and Stirling's formula.

Following the method of Theorem~\ref{lagtheo} we have
\begin{eqnarray*}
N_{2n,z}&\leq & s_\theta^{-1} \int_{-\infty}^\infty\sum_{r,s=0}^n
 |b_{n,r}b_{n,s}|s_\theta^{-2r-2s}
x^{2r+2s}\rme^{-x^2}\,\rmd x \\
&\leq & s_\theta^{-4n-1} \sum_{r,s=0}^n
 |b_{n,r}b_{n,s}|\pi^{1/2} 2^{r+s}r!\,s!\\
&\leq & s_\theta^{-4n-1}2^{-2n}(2n)! \sum_{r,s=0}^n\frac{
2^{3r+3s}r!\,s!}{(n-r)!\, (2r)!\, (n-s)!\, (2s)! }\\
&\leq & s_\theta^{-4n-1}(n+1)^{-1/2}\left(\sum_{r=0}^n\frac{
2^{3r}r!\,n!}{(n-r)!\, (2r)! }\right)^2\\
&\leq & s_\theta^{-4n-1}(n+1)^{-1/2}\left(\sum_{r=0}^n2^{3r}\frac{
n!}{(n-r)!\, r! }\, \frac{(r!)^2}{(2r)!}\right)^2\\
&\leq & s_\theta^{-4n-1}(n+1)^{-1/2}2^{2n}
\left(\sum_{r=0}^n2^{3r}\frac{(r!)^2}{(2r)!}\right)^2\\
&\leq & s_\theta^{-4n-1}(n+1)^{-1/2}2^{2n}
\left(\sum_{r=0}^n2^{r}\sqrt{\pi(r+1)}\right)^2\\
&\leq & s_\theta^{-4n-1}\pi (n+1)^{1/2} 2^{4n+2}.\\
\end{eqnarray*}

\section{The Spectral Expansion}

Let $H$ denote the non-self-adjoint harmonic oscillator acting in
$L^2(\R)$ , with eigenvalues $\lam_n=z^2(2n+1)$ and spectral
projections $P_n$. If the right hand-side of the expansion
\be%
\rme^{-Ht}=\sum_{n=0}^\infty \rme^{-\lam_n t} P_n\label{ser}
\ee%
is norm convergent, then by comparing the action of the two sides
on the eigenfunctions $\phi_n$ we see that they coincide on a
dense subspace, and hence on the whole of $L^2(\R)$.

If we put
\[
s_z=\limsup_{n\to\infty} n^{-1}\log(\norm P_n\norm )
\]
then our theorems imply that $0<s_z <\infty$ provided
$0<|\theta|<\pi/4$. They also provide explicit upper and lower
bounds on $s_z$.

\begin{theorem}
The spectral expansion (\ref{ser}) is norm convergent if
$t>t_z=s_z/(2\cos(2\theta))$ and is norm divergent if $0\leq t
<t_z$.
\end{theorem}

\Proof For $t>t_z$ the terms of the series decrease at an
exponential rate, while for $0\leq t <t_z$ they are not uniformly
bounded in norm.

\section{Numerical Results}

The non-self-adjoint harmonic oscillator (\ref{harmonic}) has
eigenvalues $\lam_n=z^2(2n+1)$ and eigenfunctions
\[
\phi_n(x)=k_n\rme^{-z^2x^2/2}H_n(zx)
\]
for $n=0,1,...$, where $k_n$ are normalization constants, $H_n$
are the Hermite polynomials, and $|\arg(z)|<\pi/4$.

\begin{theorem}\label{asymp}
If $P_n$ is the $n${\rm th} spectral projection of $H$ and
$z=r\rme^{i\theta}$ then
\[
\liminf_{n\to\infty} n^{-1}\log(\norm P_n\norm)\geq
\log(\sec(2\theta)).
\]
\end{theorem}

\Proof This follows directly from Theorem~\ref{main} upon
observing that $\norm P_n\norm =N_{n,z}$ and
$s_\theta=\cos(2\theta)^{1/2}$.

We have previously evaluated these norms numerically for
$z^4=c=\sqrt{i}$, i.e. $\theta=\pi/16$. See $\kap_n^{(1)}$ in
Table 4 of \cite{AD}. It appears from the computations there that
\[
\lim_{n\to\infty} n^{-1}\log(\norm P_n\norm)\sim 0.40
\]
which is considerably larger than the lower bound $0.079$ of
Theorem~\ref{asymp}.

We now report on a more systematic numerical investigation of the
spectral projections of (\ref{harmonic}). We evaluated
$\sig_n(\theta)=\sqrt{\norm P_n\norm\, / \,\norm P_{n-2}\norm }$
for various $n$ and $\theta $ using Maple. (This was easier than
evaluating $ \norm P_n\norm\, / \,\norm P_{n-1}\norm $ because
different algorithms are needed for even and odd $n$.) The method
used was the same as that described in \cite[sect. 4.3]{AD}. We
put $Digits:=200$, and included enough terms of the sequence
determining the eigenvector to achieve stability. For each
$\theta$ it appeared that $\sig_n(\theta)$ was an increasing
function of $n$, so the limiting value is probably larger than the
computed value. For $\theta=0$ the operator $H$ is self-adjoint,
and the projections have norm $1$. As stated earlier one must
restrict $\theta$ to the range $| \theta | <\pi/4$. The results
are shown for $n=100$ in Table 1. The second column lists the
constants $s_\theta^{-2}=\sec(2\theta)$ (rounded down) associated
with the lower bound of Theorem~\ref{main}. The fourth column
lists the constants $4s_\theta^{-2}=4\sec(2\theta)$ (rounded up)
associated with the upper bound of Theorem~\ref{upper}. The final
column lists the values of $\mu(\theta)=\exp(\tan(2\theta))$, for
reasons explained below.

\begin{table}[h]  
\[
\begin{array}{ccccc}
\theta  & s_\theta^{-2} & \sig_{100}(\theta)
 &4s_\theta^{-2} &\mu(\theta)\\
\cline{1-5}
0       &1      & 1       &  4     &1 \\
0.025   &1.012  & 1.165   & 4.050  & 1.172  \\
0.05    &1.051  & 1.369   & 4.206  & 1.384  \\
0.1     &1.236  & 1.953   & 4.945  & 2.068  \\
0.15    &1.701  & 3.062   & 6.806  & 3.961   \\
0.20    &3.236  & 6.282   & 12.945 & 21.708\\
\end{array}
\]
\begin{center}
Table 1
\end{center}
\end{table}

The approximations $\mu(\theta)$ were obtained by the following
non-rigorous method. For even values of $n$ the eigenfunction
$\phi_n$ of $H$ is an even function of $x$ which is concentrated
around the points $\pm x_0$, where $x_0$ is defined below. On the
positive half-line the semi-classical analysis of \cite[Sect.
2]{EBD} suggests that for large enough $\eta >0$
\[
\phi(s+x_0)\sim\rme^{-\psi_1s-\psi_2s^2/2}
\]
is an approximate eigenvector of $H$ with approximate eigenvalue
$\lam$, where $x_0=\eta$, $\psi_1=i\eta$, $\psi_2=-iz^4$, and
$\lam=(1+z^4)\eta^2$; in the notation of \cite{EBD} we are putting
$c=z^4$ and $\alp=1$, and are ignoring the term involving
$\psi_3$.

If $n$ is a positive integer and we put $\eta=\{
n/\cos(2\theta)\}^{1/2}$, then a direct calculation shows that
$\lam=2nz^2$, which equals the $n$th eigenvalue of $H$ to leading
order as $n\to\infty$. This suggests that
\begin{eqnarray*}
\norm P_n\norm &\sim & \frac{\int_0^\infty |\phi(s+x_0)|^2\,\rmd s
}
{|\int_0^\infty \phi(s+x_0)^2\,\rmd s| }\\
&\sim & \frac{\int_{-\infty}^\infty
\rme^{-2\Re(\psi_1)s-\Re(\psi_2)s^2}\,\rmd s }
{|\int_{-\infty}^\infty \rme^{-2\psi_1s-\psi_2s^2}\,\rmd s| }\\
&=&\exp\{ n\tan(2\theta)\}.
\end{eqnarray*}
In view of the crude character of the approximations above, the
similarity of $\sig_{100}(\theta)$ and $\mu(\theta)$ in Table 1 is
interesting. We conjecture that a more detailed semiclassical
analysis might yield the correct asymptotic constant. This also
seems the best hope for treating more general non-self-adjoint
\Schrodinger operators.

\end{document}